# COMBINATORICS OF A FRACTAL TILING FAMILY

Hassan Douzi,
Université Ibn Zohr, Faculté des sciences, BP8106
Département de mathématiques, Agadir, Morocco.
Email: douzi_h @yahoo.fr

## Abstract

In this paper, we propose to enumerate all different configurations belonging to a specific class of fractals: A binary initial tile is selected and a finite recursive tiling process is engaged to produce auto-similar binary patterns. For each initial tile choice the number of possible configurations is finite. This combinatorial problem recalls the famous Escher tiling problem [2]. By using the Burnside lemma we show that there are exactly 232 really different fractals when the initial tile is a particular 2×2 matrix. Partial results are also presented in the 3×3 case when the initial tile presents some symmetry properties.

## 1 Introduction

Very simple processes (IFS, L-Systems, Chaos ...etc.) are often used to generate complex fractal patterns, and make very rich the fractal world. In this paper we generate a family of fractals, designed by "fractal tiling", by using a very simple algorithm: Each pattern is initially determined by the choice of an initial binary square matrix tile and the association of each matrix element to a transformation chosen in the planar symmetric group. Then a super-tile is constructed by replacing each element of the matrix by its associated transformation applied on the initial tile. Finally a recursive tiling repeats the same operation and reproduces larger super-tiles with auto-similar properties.

The interest of such fractal class resides in the fact that, for each initial tile, the number of possible configurations is always finite. So, by enumerating them, we can have a precise idea about all kinds of fractal patterns composing them. Also, with this method, we can construct every fractal with auto-similar matrix patterns. This is, for example, the case of a large part of IFS (Iterated Function System) fractals which use the same planar symmetry group transformations [1]. This enumeration problem recalls the famous Escher tiling problem where the first super-tile is simply repeated periodically in the plane [2] [6]. Many publications have treated combinatorial aspects of this tiling problem especially by using the Burnside lemma [3].

Here, we use the Burnside lemma to calculate the number of different fractal tiling patterns obtained when the initial tile is a 2×2 matrix with one zero element (this is the only interesting case for this size). The number of all possible patterns is 2048 but only 232 are really different (See Appendix for an illustration of all of them). After that we consider the case where the initial tile is a 3×3 matrix. The number of possible cases explodes to: 387420489. We limit our investigations to the case where the zero distribution in the initial tile respects a diagonal symmetry (for only one diagonal). This choice is motivated, besides its aesthetic symmetrical aspect, by the fact that it is a generalization of the 2×2 case. We calculate the number of different fractal patterns in the same manner, by using the Burnside Lemma, and we get exactly: 11043660.

## 2. A fractal tiling method

Firstly we consider a binary square matrix $M_0$ with size n×n. we associate, to every nonzero element a transformation selected in the planar symmetry group G. This group is composed of transformations which act on plane patterns by rotation (with successive quarter turns) or mirror reflection. G is completely determined by the eight following transformations [4] [5]:

- R0=Id, R1, R2 and R3: respectively rotations, counterclockwise, with angles 0 °, 90 °, 180 ° and 270 °.
- K0, K1, K2 and K3: respectively the horizontal mirror reflections of previous rotations, i.e.: M o R0, M o R1, M o R2 and M o R3, where M is the horizontal mirror transformation.

Then we build a super-tile $M_1$, composed of a binary square matrix of size $n^2 \times n^2$, as follows:
- Every zero element of $M_0$ is replaced by a null matrix of size $M_0$.
- Every nonzero element of $M_0$, associated to a transformation $T \in G$, is replaced by the matrix $T(M_0)$.

Finally we end the construction by a recursive tiling process which assembles bigger super-tiles in the same manner:
- We repeat, for a finite number of iterations, the following instructions:
  - We have a current square matrix $M_k$ of size $n^{k+1} \times n^{k+1}$ (Initially k=1)
  - Every zero element of $M_0$ is replaced by a null matrix of size $M_k$
  - Every nonzero element of $M_k$ associated to a transformation $T \in G$ is replaced by the matrix $T(M_k)$.
  - The new matrix of size $n^{k+2} \times n^{k+2}$ is designed by $M_{(k+1)}$

At the end we get a binary matrix with auto-similar patterns. By continuing the process to infinity, we can identify them to authentic fractals. For example to generate the famous Sierpinsky triangle we consider $M_0 = \begin{pmatrix} 0 & 1 \\ 1 & 1 \end{pmatrix}$ and we associate all its nonzero elements to the identity transformation (Id=R0) (*Figure 1*).

We note that when all elements in $M_0$ are nonzero we get monocolor matrix $M_k$ without any discernable fractal pattern. But if we associate to each matrix element a value corresponding to its associated transformation we get an image with fractal textural patterns *(Figure 1)*. But we won't consider this aspect now in this study.

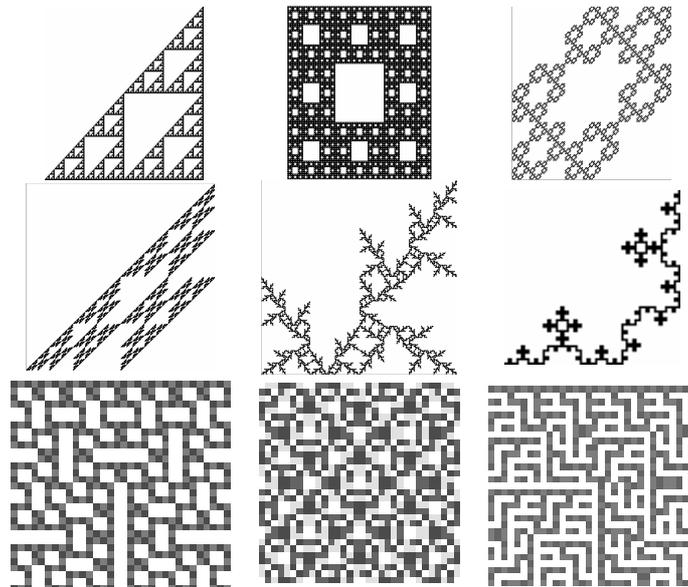

***Figure 1*** *: Panorama of fractal tiling : Upward three classical fractals: Sierpinsky triangle and square with respectively $M_0$=(0 Id ;Id Id) and $M_0$=(Id Id Id ;Id 0 Id ; Id Id Id) , Von-Koch flakes with $M_0$=(0 Id Id ;Id 0 Id ; Id Id 0). In the middle other less known fractals with (from left to right) $M_0$=(0 Id ; Id R2) (maple leaf), $M_0$=(0 Id ; R3 K2), $M_0$=(0 0 Id ;0 0 R1 ; Id R3 R2 ). Downward textural fractals generated when all element in $M_0$ are nonzero (from left to right): $M_0$=( Id Id ; K3 K1 ) $M_0$=( Id R1 ; K2 K3 ) $M_0$=( Id Id ; K2 K1 ) binary matrix are replaced by multi-values matrix associated to symmetric group transformations.*

# 3. The 2×2 fractal tiling

When $M_0$ is a (2×2) matrix we have only five different zeros distributions (up to a rotation):

$$\begin{pmatrix} 0 & 0 \\ 0 & 0 \end{pmatrix} \begin{pmatrix} 0 & 0 \\ 0 & 1 \end{pmatrix} \begin{pmatrix} 0 & 0 \\ 1 & 1 \end{pmatrix} \begin{pmatrix} 0 & 1 \\ 1 & 1 \end{pmatrix} \begin{pmatrix} 1 & 1 \\ 1 & 1 \end{pmatrix}$$

The first three distributions (from the left) are obviously not interesting because thy lead either to empty patterns (with only zeros) or some dusty patterns which recalls Cantor sets. The last distribution, where all elements are nonzero, generates uniform matrixes (with only ones) but if we assign different values to transformations associated to $M_0$ elements we get fractal textures (*figure 1*). So we will focus our study on the fourth configuration with only one zero element.

Each nonzero element in this case may be associated with one of the 8 possible transformations in G and the zero has 4 possible positions in $M_0$. Therefore the total number of possible configurations in this case is: $8^3*4 = 2048$. We will denoted every configuration by a quadruplet (a,b,c,d) which verify:

- Only one element in the quadruplet is zero.
- "a" is either the transformation associated to $M_0(1,1)$ or zero.
- "b" is either the transformation associated to $M_0(1,2)$ or zero.
- "c" is either the transformation associated to $M_0(2,1)$ or zero.
- "d" is either the transformation associated to $M_0(2,2)$ or zero.

In the next sections we will consider that all configurations obtained by plane rotation or mirror reflection, from the same initial one, are identical and we propose to enumerate all different configurations in this sense.

## *3-1. Redundancy elimination*

When $M_0 = \begin{pmatrix} 0 & 1 \\ 1 & 1 \end{pmatrix}$ (up to quarter turns rotation) the distribution of nonzero elements presents a diagonal symmetry which involves the existence of many redundant configurations. So we must begin by enumerating repeated configurations and eliminate their recurrence in the set of possible configurations.

Let consider configurations associated to the quadruplets which have the form (0, b, c, d) (we don't lose in generality because we can get all the others configurations simply by rotation). An exhaustive browsing of the quadruplet set with this form, with a computer program, gives us eight classes containing each of them eight redundant configurations:

- **(0 , R0 or K1 , R0 or K1 , R0 or K1)** : (0,R0,R0,R0 ) , (0,R0,K1,R0) , (0,K1,R0,R0) , (0,K1, K1,R0) , (0,R0,R0,K1 ) , (0,R0,K1,K1) , (0,K1,R0,K1) , (0,K1, K1,K1)
- **(0 , R0 or K1 , R0 or K1 , R2 or K3)** : (0,R0,R0,R2 ) , (0,R0,K1,R2) , (0,K1,R0,R2) , (0,K1, K1,R2) , (0,R0,R0,K3 ) , (0,R0,K1,K3) , (0,K1,R0,K3) , (0,K1, K1,K3)
- **(0 , R2 or K3 , R2 or K3 , R0 or K1)** : (0,R2,R2,R0 ) , (0,R2,K3,R0) , (0,K3,R2,R0) , (0,K3, K3,R0) , (0,R2,R2,K1 ) , (0,R2,K3,K1) , (0,K3,R2,K1) , (0,K3, K3,K1)
- **(0 , R2 or K3 , R2 or K3 , R2 or K3)** : (0,R2,R2,R2 ) , (0,R2,K3,R2) , (0,K3,R2,R2) , (0,K3, K3,R2) , (0,R2,R2,K3 ) , (0,R2,K3,K3) , (0,K3,R2,K3) , (0,K3, K3,K3)
- **(0 , R1 or K0 , R3 or K2 , R0 or K2) :** (0,R1,R3,R0 ) , (0,R1,K2,R0) , (0,K0,R3,R0) , (0,K0, K2,R0) , (0,R1,R3,K1 ) , (0,R1,K2,K1) , (0,K0,R3,K1) , (0,K0, K2,K1)
- **(0 , R1 or K0 , R3 or K2 , R2 or K3) :** (0,R1,R3,R2 ) , (0,R1,K2,R2) , (0,K0,R3,R2) , (0,K0, K2,R2) , (0,R1,R3,K3 ) , (0,R1,K2,K3) , (0,K0,R3,K3) , (0,K0, K2,K3)
- **(0 , R3 or K2 , R1 or K0 , R0 or K1) :** (0,R3,R1,R0 ) , (0,R2,R1,R0) , (0,R3,K0,R0) , (0,K2, K0,R0) , (0,R3,R1,K1 ) , (0,K2,R1,K1) , (0,R3,K0,K1) , (0,K2, K0,K1)

- **(0 , R3 or K2 , R1 or K0 , R2 or K3) :** (0,R3,R1,R2 ) , (0,K2,R1,R2) , (0,R3,K0,R2) , (0,K2, K0,R2) , (0,R3,R1,K3 ) , (0,K2,R1,K3) , (0,R3,K0,K3) , (0,K2, K0,K3)

We notice that all those redundant configurations are, in fact, characterized by the following properties:

*Corollary 1*
*Redundant configurations of the form (0, b, c, d) are characterized by the three next conditions:*
1) *For matrix $M_k$ ($k≥0$), associated to Redundant configurations, the action of all G transformations are two by two equivalent: $R0(M_k)=K1(M_k)$ ; $R1(M_k)=K0(M_k)$ ; $R2(M_k)=K3(M_k)$ and $R3(M_k)=K2(M_k)$. In other word, for redundant configurations, there are only four different transformations in the symmetry group G. This gives us eight configurations for each redundancy because we can change each transformation by its equivalent one.*
2) *Non diagonal elements (($M_0(1,2)$ and $M_0(2,1)$)) are associated with transformations which are symmetrical each other with respect to the diagonal, in other words they are limited to the next 16 possibilities:*
$$(( R0 \text{ or } K1 ) \text{ and } ( R0 \text{ or } K1 )) ; (( R2 \text{ or } K3 ) \text{ and } (R2 \text{ or } K3 )) ;$$
$$(( R1 \text{ or } K0 ) \text{ and } ( R3 \text{ or } K2 )) ; ((R3 \text{ or } K2 ) \text{ and } (R1 \text{ or } K0 ))$$
3) *The diagonal element $M_0(2,2)$ is only associated to transformations (R0 or K1) or (R2 or K3) i.e. to transformations which are their own symmetrical for the condition 2.*

A more formal proof of these properties will be given in the next section (*Proposition 2*).
The same arguments to identify redundant configurations stay true when we change the zero position in $M_0$. Then we have a total number of 8*4=32 classes with each one containing 8 identical configurations. If we take only one representative for each class the configuration set without redundancy, noted by $E_2$, contains: 2048 - 4*64 + 32 = 1824 configurations.

## 3-2. Enumeration of different configurations

The second step, in the enumeration of different configurations, consists in identifying all configurations obtained from the same one by rotation or by mirror reflection (or a combination of both of them). More formally we must enumerate all the equivalence classes associated with the group action of G on the set $E_2$. The Burnside Lemma (called also Frobenius-Cauchy Lemma or Polya Theorem) is particularly useful for accomplishing this task [3]:
For every element $x \in E_2$, we define the x-Orbit as the set $O_x = \{T(x)/T \in G\}$. The set of all Orbits, denoted by (G / E2), forms a partition of $E_2$ which can be associated with an equivalence relation ($x \sim y \Leftrightarrow y \in O_x$), (G / E2) is then the quotient set. The set of points fixed by a transformation $T \in G$ is denoted by $E_T = \{x \in E_2 / T(x) = x\}$. The Burnside Lemma can then be formulated as a relation between the number of $E_2$-orbits and the number of all points fixed by a G transformation:

$$|E_2 / G| = \frac{\sum_{T \in G} |E_T|}{|G|}$$

*Proposition1:*
*The number of different configurations for the group action of G on the configuration set $E_2$ is 232.*

*Proof*
Among the eight transformations of the symmetry group G, only the sets $E_{R0}$, $E_{K1}$ and $E_{K3}$ are nonempty. Indeed the existence of one zero element in $M_0$ imply the existence of one null quadrant in every matrix $M_k$. We note then that transformations R1, R2, R3, K0 and K2 displace, in their action, the null quadrant position of $M_k$ while R0, K1 and K3 let it in the same diagonal position.

The transformation R0=Id let, of course, invariant all $E_2$ elements. For K1 and K3 we note that invariant configurations belong to redundant configuration classes studied before. Indeed the only invariant configurations for K1 (respectively K3) are those which are symmetrical (in the sense of condition 2 corollary 1) with respect to the principal diagonal (respectively the second diagonal). A more general proof will be presented in *Proposition3)*. Then K1 (respectively K3) let invariant 16 element of $E_2$. Finally the Burnside formula can be applied as follows:

$$|E_2/G| = \frac{|E_{R0}|+|E_{K1}|+|E_{K3}|}{|G|} = \frac{1824+16+16}{8} = 232 \blacksquare$$

A representative list of all those 232 different configurations is presented in the *Appendix*.

# 4. The 3×3 fractal tiling case

Now we consider $M_0$ as a binary square 3×3 matrix where every nonzero element is associated to a transformation belonging to the group G. The total number of possible configurations is given by the next sum: $\sum_{k=0}^{k=9} \binom{k}{9} 8^k = 9^9 = 387420489$; where each term of the sum represents the number of possible configurations for a given number of nonzero elements in $M_0$.

We have seen, in the 2×2 case, that the study of the zeros distribution and their symmetrical disposition in $M_0$ are essential to identify similar configurations (up to a rotation or mirror reflection). In the 3×3 case the difficulty comes from the high number of distributions to be handled: to dispose k zeros in $M_0$ we have $\binom{k}{9}$ possible choices, so we have $\sum_{k=0}^{k=9} \binom{k}{9} = 2^9 = 512$ possible distribution.

## *4-1. Configurations with diagonal symmetry*

We will limit our investigations to configurations which present, in $M_0$, a symmetrical zero distribution with respect the one of the two $M_0$ diagonals and in the same time an asymmetrical zero distribution with respect to the other diagonal. We denote by *DS* the set of this class of configurations. This limitation is justified, beside the esthetical aspect of symmetry, by the fact that it reproduce and generalize the situation studied in the 2×2 case. We can then use, in the same manner, the Burnside Lemma. For those special configurations we have two separate situations:
1) The zero distribution is only symmetrical with respect to the principal diagonal (diagonal from up-left to down-right) and this subfamily is denoted by *DS1*
2) The zero distribution is only symmetrical with respect to the second diagonal and this subfamily is denoted by *DS2*

We can switch between the two situations by a simple 90° rotation so we can only study one of them. Let *N1* be the set of configurations which present symmetrical zero distributions with respect to the principal diagonal (we have $DS1 \subset N1$). Then the number of configurations in *N1* depends on the number of zero distributions in the principal diagonal of $M_0$ (among 3 elements) and the number of them situated upward the diagonal (also among 3 elements), the symmetry determine the position of the remaining zeros. Knowing that every nonzero element is associated to a transformation belonging to the group G, the *N1* cardinal *is given by:*

$$|N1| = \left(\sum_{k=0}^{k=3}\binom{k}{3}8^k\right)\left(\sum_{l=0}^{l=3}\binom{l}{3}8^{2*l}\right) = 9^3 * 65^3 = 200201625$$

*DS1* is a subset of *N1* and we have: *N1=DS1 U DS12* where *DS12* denotes the set of configurations which present in $M_0$ a symmetrical zero distribution with respect to the two diagonals. *DS12* contains 16 possible zeros distribution because we have 4 possibilities of symmetrical zeros distribution for diagonal elements and the same number of possibilities for non diagonal elements:

$$\begin{pmatrix}0&0&0\\0&0&0\\0&0&0\end{pmatrix}\begin{pmatrix}0&0&1\\0&0&0\\1&0&0\end{pmatrix}\begin{pmatrix}0&1&0\\1&0&1\\0&1&0\end{pmatrix}\begin{pmatrix}0&1&1\\1&0&1\\1&1&0\end{pmatrix}\begin{pmatrix}0&0&0\\0&1&0\\0&0&0\end{pmatrix}\begin{pmatrix}0&0&1\\0&1&0\\1&0&0\end{pmatrix}\begin{pmatrix}0&1&0\\1&1&1\\0&1&0\end{pmatrix}\begin{pmatrix}0&1&1\\1&1&1\\1&1&0\end{pmatrix}$$

$$\begin{pmatrix}1&0&0\\0&0&0\\0&0&1\end{pmatrix}\begin{pmatrix}1&0&1\\0&0&0\\1&0&1\end{pmatrix}\begin{pmatrix}1&1&0\\1&0&1\\0&1&1\end{pmatrix}\begin{pmatrix}1&1&1\\1&0&1\\1&1&1\end{pmatrix}\begin{pmatrix}1&0&0\\0&1&0\\0&0&1\end{pmatrix}\begin{pmatrix}1&0&1\\0&1&0\\1&0&1\end{pmatrix}\begin{pmatrix}1&1&0\\1&1&1\\0&1&1\end{pmatrix}\begin{pmatrix}1&1&1\\1&1&1\\1&1&1\end{pmatrix}$$

The number of *DS12* configurations is then given by the next formula:

$$|DS12| = \left(\sum_{k=0}^{k=3} 8^k\right)\left(\sum_{l=0}^{l=3} 8^{2*l}\right) = \frac{8^4-1}{8-1} * \frac{64^4-1}{64-1} = 155788425$$

*DS1* and *DS12* are disjoint so the cardinal *of DS1* (which is also the cardinal of *DS2*) is given by:

$$|DS1| = |DS2| = |N1| - |DS12| = \sum_{k=0}^{k=3}\sum_{l=0}^{l=3}\left(\binom{k}{3}\binom{l}{3}-1\right)*8^{(k+2*l)} = 44413200$$

Where every term indexed by *(k,l)* represents the number of configurations with a distribution of *(k+2\*l)* nonzero elements in $M_0$. Finally the number of configurations in *DS* is given by:

$$|DS| = |DS1| + |DS2| = 88826400$$

## *4-2. Redundancy elimination*

Like the 2×2 case, the set *DS* contains many redundant configurations due to the symmetry of zero distribution in $M_0$. In *DS1* (and equivalently in *DS2*) the redundancy is characterized by the next result.

*Proposition2*
*Redundant configurations in DS1 are characterized by the three following conditions:*
   *(1) G transformations are two by two equivalent for $M_k$ (k≥0) matrix:*
       *R0($M_k$)=K1($M_k$) ; R1($M_k$)=K0($M_k$) ; R2($M_k$)=K3($M_k$) and R3($M_k$)=K2($M_k$)*
   *(2) Diagonal elements ($M_0(i,j)$ with 1≤i=j≤3) are associated uniquely to the four next transformations: R0 ≡K1 or R2 ≡K3.*
   *(3) If $M_0(i,j)$ with 1≤i<j≤3 is associated to the transformation (R0 or K1) (respectively (R1 or K0), (R2 or K3), (R3 or K2)) then $M_0(j,i)$ is associated to the transformation (R0 or K1) (respectively (R3 or K2), (R2 or K3), (R1 or K0)).*

*Proof*
Firstly a configuration is redundant if and only if we can substitute at least one transformation, associated to a given element of $M_0$, without any changing in the obtained configuration. In other words $M_0$ and by recurrence $M_k$ (k≥0) give the same results under the action of, at least, two different transformations in the group G. The *DS1* configurations are characterized by a symmetry of zeros distribution only with respect to the principal diagonal so (exactly like in the case 2×2). The condition (1) gives the list of transformation couples susceptible to give the same result on a

*DS1* configuration. For example we can verify that for R0 only K1 can give the same result on a *DS1* configuration because it is the only one to preserve the same zero distribution (and the same symmetry direction)....etc.

Secondly, in a given *DS1* configuration, if $M_0(i,i)$ (for i=1,2 or 3) is associated to a transformation T which change the direction of the symmetry diagonal in $M_0$ (i.e. associated to T= R1, R3, K0 or K2) then the matrix $M_1$ wont presents anymore any diagonal symmetry in its zero distribution. Now if we replace T by an equivalent transformation (in the sense of condition (1)) then we will get the same matrix $M_1$ but a different matrix $M_2$ due to the loose of symmetry in $M_1$. Thus the two configurations are not identical and this establishes the condition (2).

Finally, and in the same manner, if the condition (3) is not respected in $M_0$ we will have a matrix $M_1$ without any diagonal symmetry zero distribution, so when we replace a transformation associated to $M_0$ by its equivalent one (in the sense of condition (1)) we get the same matrix $M_1$ but not the same matrix $M_2$ and this implies that the considered configuration can't be redundant ■.

The family of redundant configurations in *DS1* (respectively in *DS2)* is denoted by *DSS1* (respectively *DSS2)*. We note that *DSS2* can be obtained by rotating by 90° all *DSS1* configurations. The set of all redundant configurations in *DS* is denoted by *DSS* and we have *DSS =DSS1 U DSS2*.

From *Proposition2*, we can conclude that for a redundant configuration in *DSS1* and for every nonzero diagonal element we have the choice between 4 possible transformations to be associated with. For every upward non diagonal and nonzero element we have 8 possible choices and 2 choices for its symmetrical downward non diagonal element. Then the cardinal of *DSS1* is given by:

$$|DSS1| = \sum_{k=0}^{k=3}\sum_{l=0}^{l=3}\left(\binom{k}{3}\binom{l}{3}-1\right)*4^k*8^l*2^l = 5^3*17^3 - \left(\frac{4^4-1}{4-1}\right)\left(\frac{16^4-1}{16-1}\right) = 242760$$

Every configuration in *DSS1*, with *k* nonzero diagonal element and *2*l* nonzero non diagonal element, belongs to a class of $2^{(k+2*l)}$ identical configurations because every transformation associated to $M_0$ can be replaced by its equivalent one (condition (1)). Then the number of classes in *DSS1* is given by:

$$N2 = \sum_{k=0}^{k=3}\sum_{l=0}^{l=3}\left(\binom{k}{3}\binom{l}{3}-1\right)*\frac{4^k*8^l*2^l}{2^{k+2l}} = 3^3*5^3 - \left(\frac{2^4-1}{1}\right)\left(\frac{4^4-1}{3}\right) = 2100$$

In the same manner, the same numbers are obtained for the *DSS2* family. The cardinal of *DSS* is then given by $|DSS| = |DSS1| + |DSS2| = 485520$ and the number of identical configuration classes is $2*N2 = 4200$. Finally by eliminating redundancy from *DS* we get a set of configurations denoted by *DSR* and its cardinal is given by : $|DSR| = |DS| - |DSS| + 2*N2 = 88345080$.

### *4-3. Enumeration of different configurations*

The second step consists, like in the 2×2 case, to identify in *DSR* all configuration pairs linked each other by a transformation of the group G. In other words enumerate all equivalence classes associated with the group action of G on the set DSR. To do it, by using the Burnside Lemma, we must calculate the fixed points for every transformation in G:

*Proposition 3*
*The number of different configurations for the group action of G on the configuration set DSR is:*
11043660.

*Proof*

As in the 2×2 case, the only transformations susceptible to have fixed points in *DSR* are R0, K1 and K3. Indeed *DSR* configurations are characterized by a symmetrical zero distribution with respect to one diagonal and an asymmetrical distribution with respect the other diagonal and the action of the transformations R1, R2, R3, K0 and K2 changes the direction of the diagonal symmetry while the action of R0, K1 and K3 preserves the same symmetrical diagonal.

The transformation R0=Id let, of course, invariant all *DSR* configurations. For transformations K1 and K3 we note firstly, according to (proposition2, condition 1), that redundant configurations in *DS1* are invariant under the action of K1 and, equivalently, those of *DS2* are invariant under the action of K3. Secondly the only configurations susceptible to be invariant are those belonging to redundant configurations classes: indeed if not (exactly as explained in proposition2) the matrix $M_1$ will lose the diagonal symmetry of zero distribution (the symmetry is however present in $M_0$) and this imply that : K1($M_k$)≠$M_k$ and K3($M_k$)≠$M_k$ for k≥2. Finally K1 and K3 let invariant, each of them, 2100 element of *DSR* because every redundant configuration class in *DS* is represented by one configuration in *DSR*.

We can now apply the Burnside formula as follows:

$$|DSR/G| = \frac{|DSR_{R0}| + |DSR_{K1}| + |DSR_{K3}|}{|G|} = \frac{88345080 + 2100 + 2100}{8} = 11043660$$

We use the same P*roposition1* notations: $DSR_{R0}$, $DSR_{K1}$, $DSR_{K3}$ are the fixed point sets of R0, K1 et K3 and DSR/G the quotient set of the G action on *DSR* ∎.

At last we observe that in fact the *proposition3*, independently from the size of $M_0$, stay valid if we have the diagonal symmetrical zero distribution as stated in the 3×3 case. Therfore, similar results can also be established for the general case n×n.

## Appendix

An exhaustive list of the 232 different ways to construct a fractal tiling, starting from an initial binary 2×2 tile with one zero element. All representative configurations have the quadruplet form (0,b,c,d) (denoted below by: motifbcd).

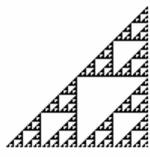 motif111
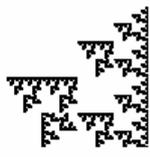 motif112
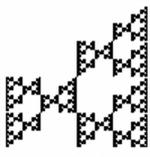 motif113
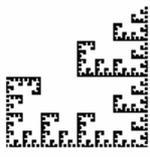 motif114
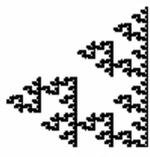 motif115
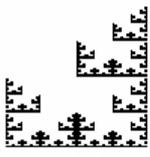 motif117
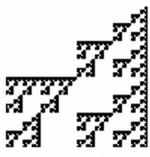 motif118

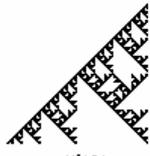 motif121
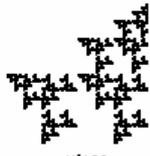 motif122
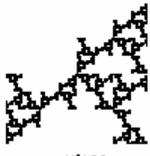 motif123
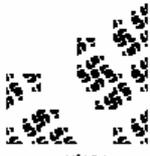 motif124
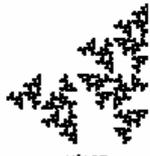 motif125
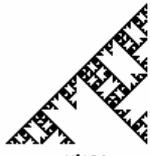 motif126
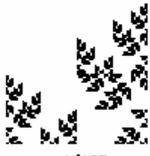 motif127

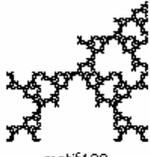 motif128
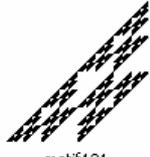 motif131
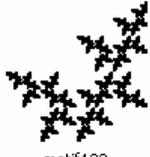 motif132
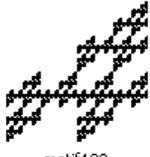 motif133
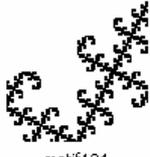 motif134
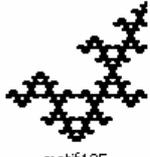 motif135
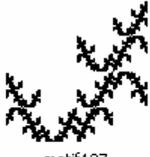 motif137

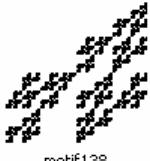 motif138
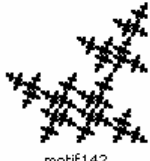 motif142
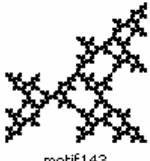 motif143
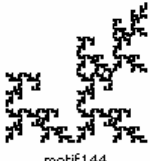 motif144
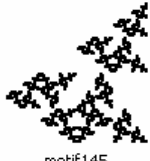 motif145
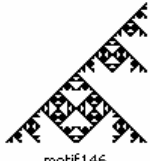 motif146
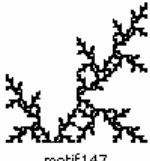 motif147

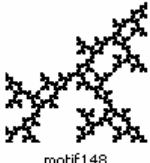 motif148
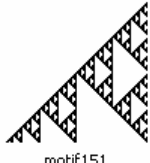 motif151
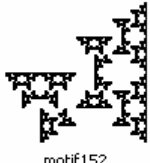 motif152
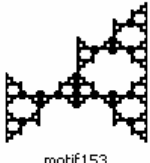 motif153
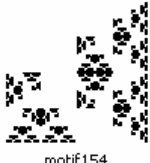 motif154
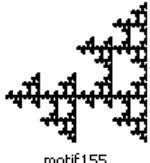 motif155
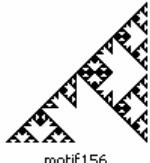 motif156

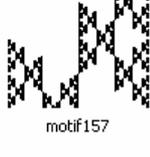 motif157
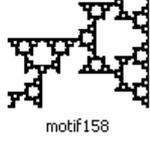 motif158
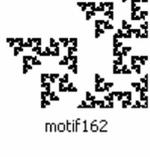 motif162
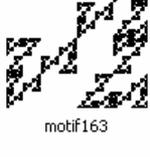 motif163
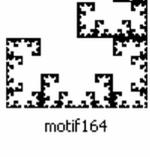 motif164
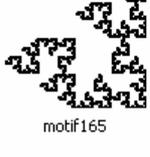 motif165
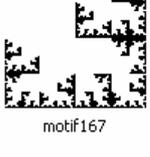 motif167

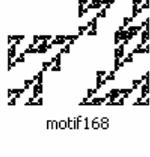 motif168
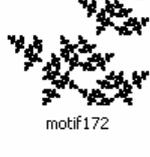 motif172
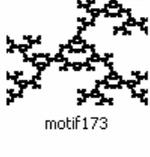 motif173
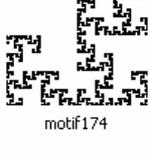 motif174
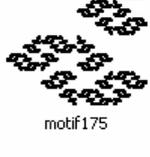 motif175
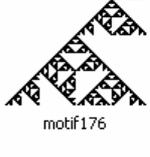 motif176
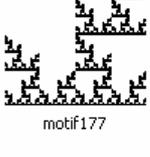 motif177

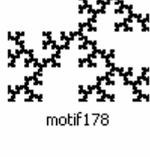 motif178
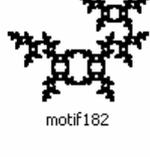 motif182
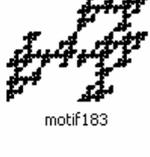 motif183
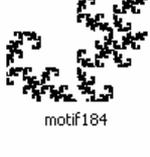 motif184
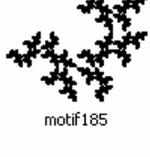 motif185
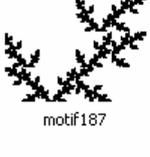 motif187
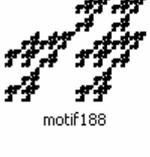 motif188

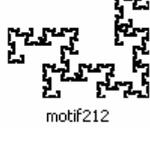 motif212
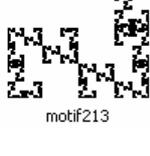 motif213
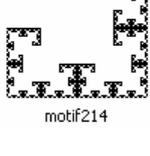 motif214
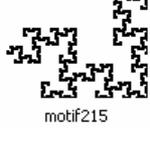 motif215
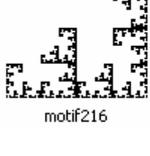 motif216
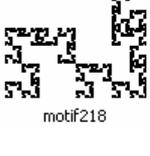 motif218
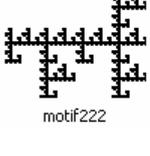 motif222

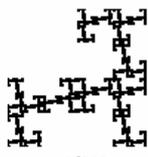 motif223
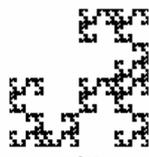 motif224
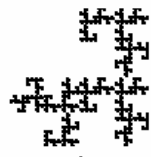 motif225
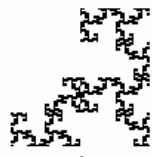 motif226
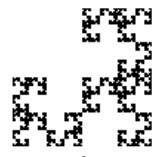 motif227
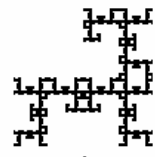 motif228
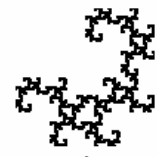 motif232

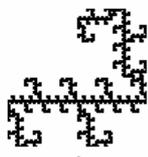 motif233
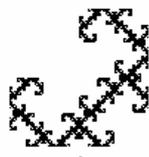 motif234
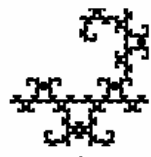 motif235
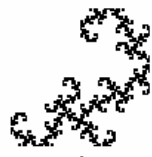 motif236
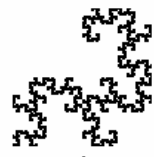 motif238
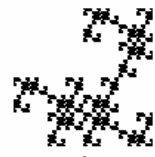 motif242
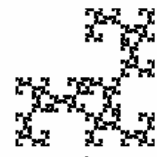 motif243

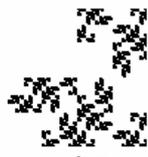 motif245
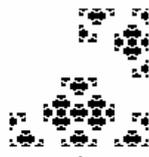 motif246
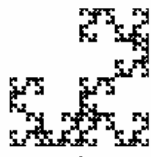 motif247
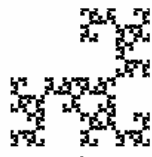 motif248
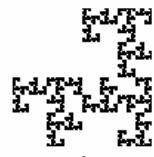 motif252
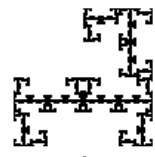 motif253
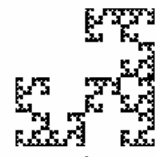 motif254

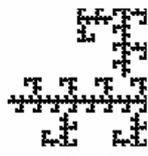 motif255
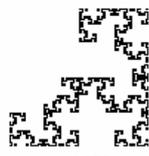 motif256
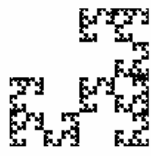 motif257
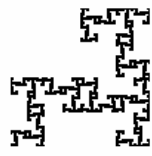 motif258
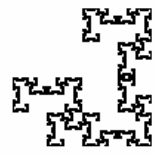 motif262
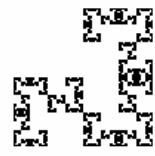 motif263
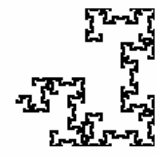 motif265

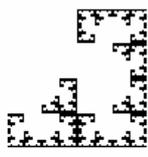 motif266
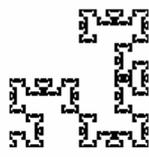 motif268
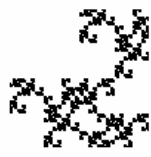 motif272
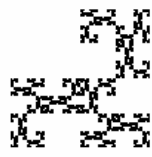 motif273
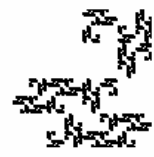 motif275
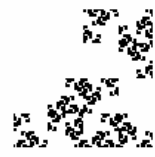 motif276
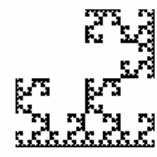 motif277

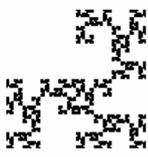 motif278
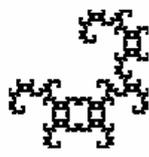 motif282
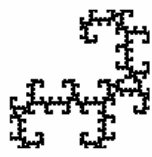 motif283
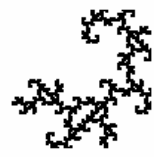 motif285
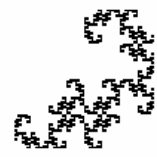 motif286
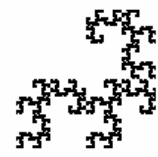 motif288
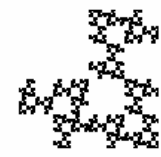 motif312

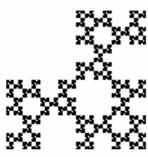 motif313
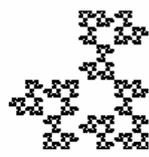 motif315
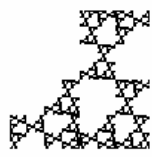 motif316
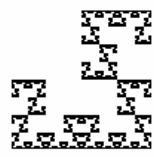 motif317
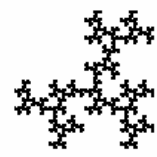 motif322
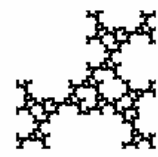 motif323
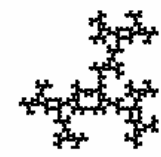 motif325

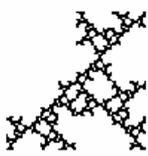 motif326
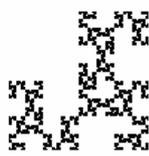 motif327
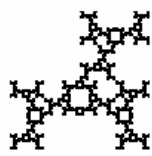 motif328
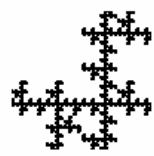 motif332
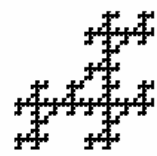 motif333
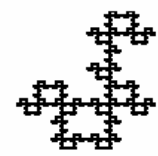 motif335
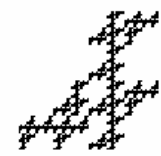 motif336

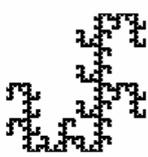 motif337
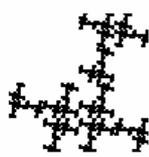 motif342
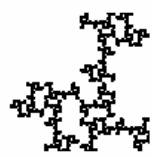 motif345
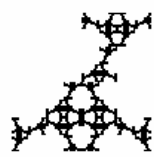 motif346
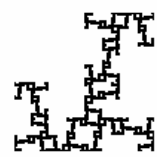 motif347
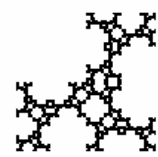 motif348
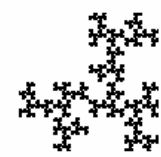 motif352

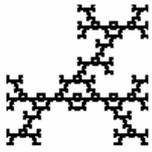 motif353
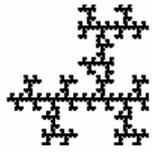 motif355
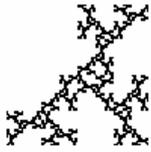 motif356
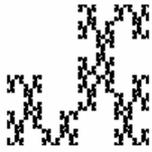 motif357
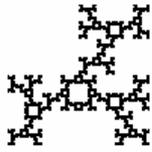 motif358
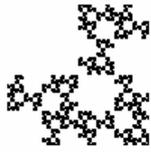 motif362
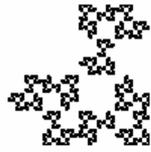 motif365
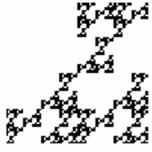 motif366
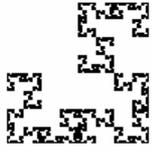 motif367
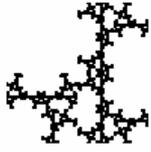 motif372
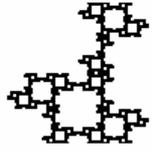 motif375
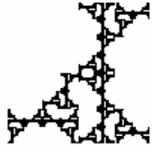 motif376
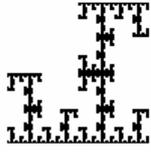 motif377
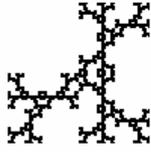 motif378
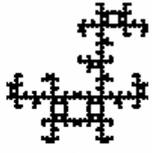 motif382
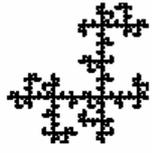 motif385
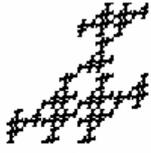 motif386
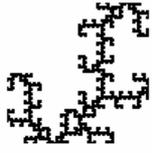 motif387
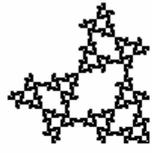 motif412
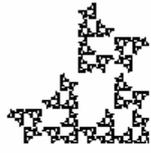 motif416
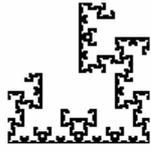 motif417
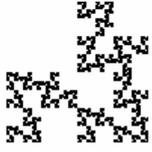 motif418
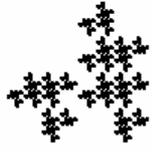 motif422
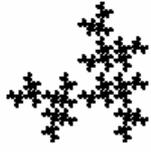 motif425
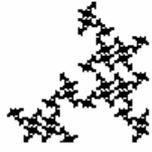 motif426
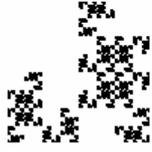 motif427
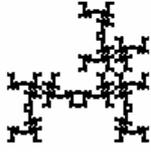 motif428
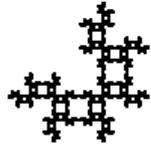 motif432
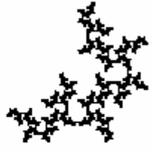 motif436
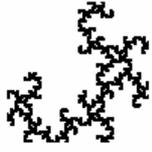 motif437
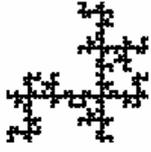 motif438
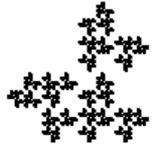 motif445
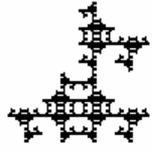 motif446
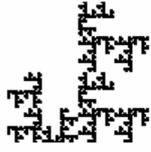 motif447
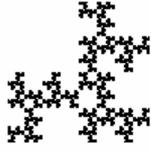 motif448
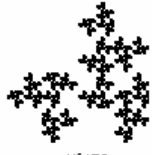 motif452
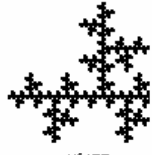 motif455
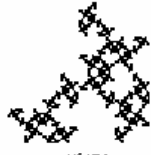 motif456
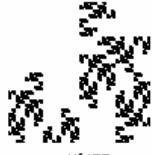 motif457
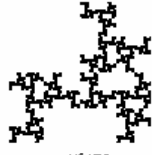 motif458
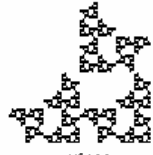 motif466
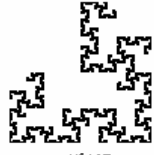 motif467
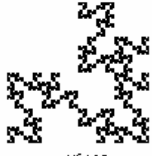 motif468
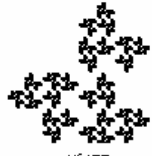 motif475
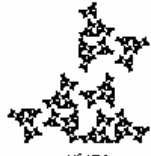 motif476
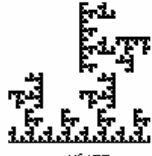 motif477
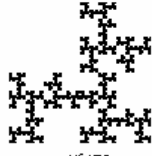 motif478
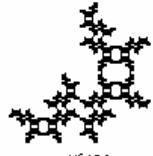 motif486
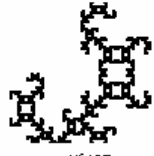 motif487
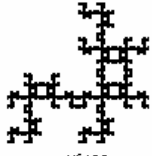 motif488
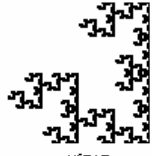 motif515
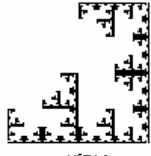 motif516
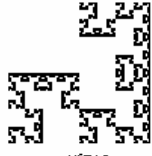 motif518
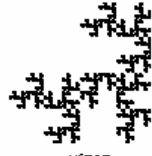 motif525
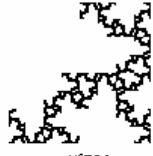 motif526
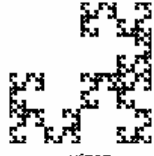 motif527
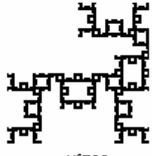 motif528
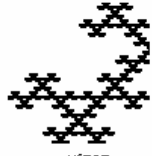 motif535
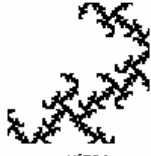 motif536
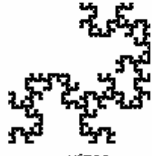 motif538
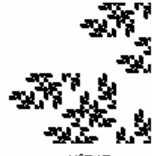 motif545
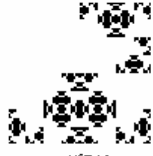 motif546
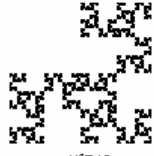 motif548

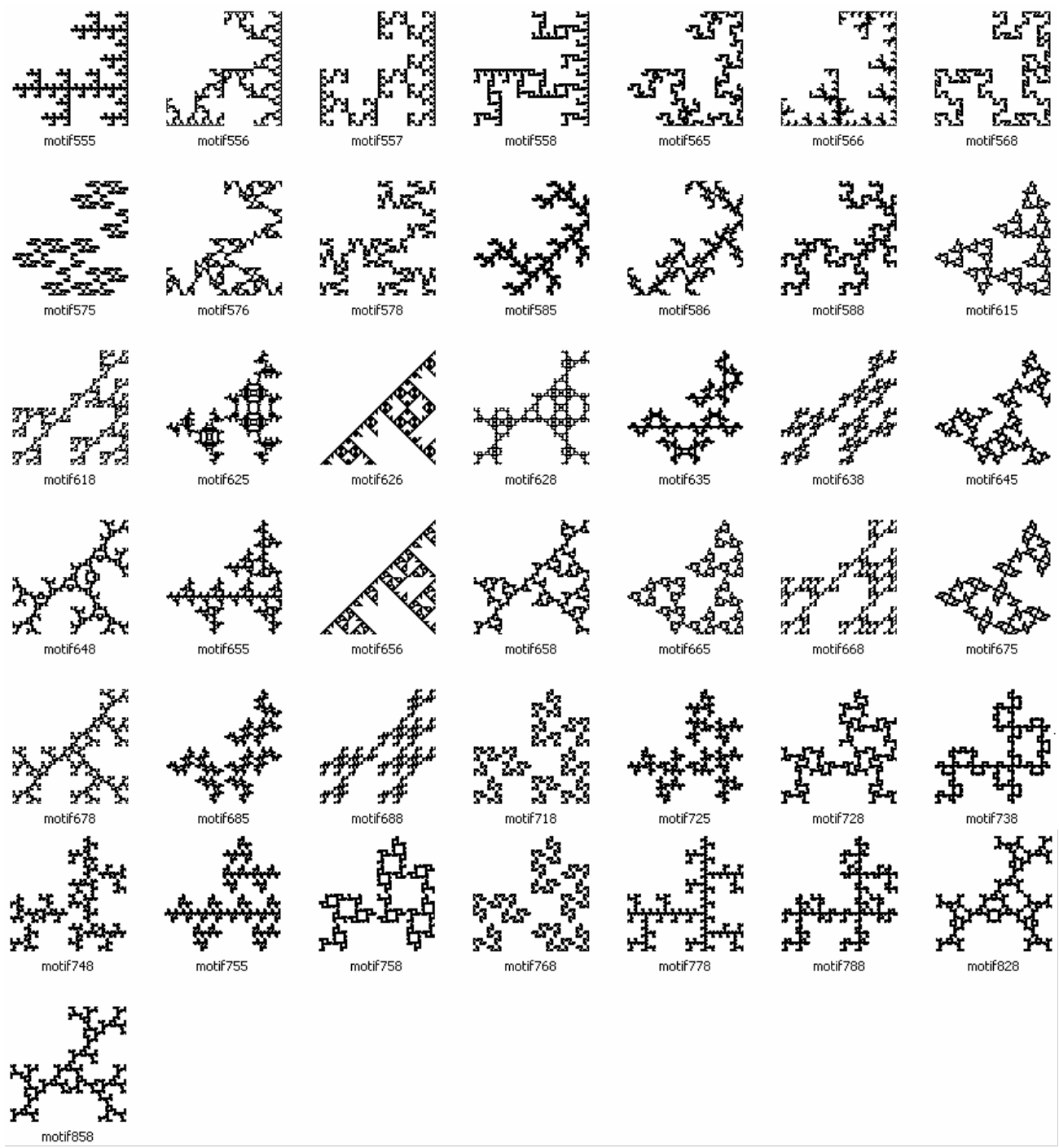

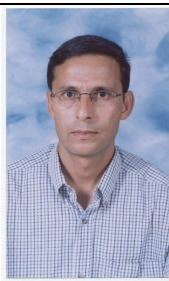

Hassan Douzi was born in Rabat, Morocco. He received the French PhD in applied mathematics from the University of Paris IX (Dauphine) in 1992. He is now a "Professeur Habilité" (research/teaching) at Ibn Zohr University in Agadir, Morocco. He has mainly worked on wavelets analysis and applications on image processing. His current research activities include also some recent fields in experimental and combinatorial mathematics.